\documentstyle[12pt]{article}
\date{May 29,  1995}

\begin{titlepage}

\title{Fuchsian groups of the second kind and representations carried by the
limit set\\[1cm]
{\small\em Dedicated to Rolf Sulanke on the occasion of his 65'th
birthday.\\[1cm]}}

\author{Ulrich Bunke\thanks{Humboldt-Universit\"at zu Berlin, Institut f\"ur
Reine Mathematik (SFB288), Ziegelstr. 13a, Berlin 10099.
E-mail:ubunke@mathematik.hu-berlin.de
} and
Martin
Olbrich\thanks{Humboldt-Universit\"at zu Berlin, Institut f\"ur Reine
Mathematik (SFB288), Ziegelstr. 13a, Berlin 10099.
E-mail:olbrich@mathematik.hu-berlin.de  }
}

\end{titlepage}

\begin{document}

\newcommand{\Bbb}{\rm}
\newcommand{\oh }{{\bf h}}
\newcommand{\om}{{\bf m}}
\newcommand{\kaaa}{{\bf k}}
\newcommand{\paaa}{{\bf p}}
\newcommand{\taaa}{{\bf t}}
\newcommand{\haaa}{{\bf h}}
\newcommand{\R}{{\bf R}}
\newcommand{\Z}{{\bf Z}}
\newcommand{\C}{{\bf C}}
 \newcommand{\Naaa}{{\bf N}}
\newcommand{\gaaa}{{\bf g}}
\newcommand{\aaaa}{{\bf a}}
\newcommand{\naaa}{{\bf n}}
 \newcommand{\db}{{\bar{\partial}}}
 \newcommand{\g}{{\gaaa}}
 \newcommand{\cZ}{{\cal Z}}
\newcommand{\cE}{{\cal E}}
\newcommand{\cC}{{\cal C}}
\newcommand{\cEp}{{{\cal E}^\prime}}
\newcommand{\cU}{{\cal U}}
\newcommand{\Hom}{{\mbox{Hom}}}
\newcommand{\End}{{\mbox{ End}}}
\newcommand{\rk}{{\mbox{rank}}}
\newcommand{\im}{{\rm im }}
\newcommand{\sign}{{\rm sign}}
\newcommand{\spann}{{\mbox{span}}}
\newcommand{\symm}{{\mbox{symm}}}
\newcommand{\cF}{{\cal F}}
\newcommand{\Ree}{{\rm Re }}
\newcommand{\Res}{{\mbox{Res}}}
\newcommand{\Imm}{{\rm Im}}
\newcommand{\inter}{{\rm int}}
\newcommand{\clo}{{\rm clo}}
\newcommand{\tg}{{\rm tan}}
\newcommand{\ee}{{\rm e}}
\newcommand{\op}{{\mbox{Op}}}
\newcommand{\tr}{{\mbox{tr}}}
\newcommand{\cs}{{c_\sigma}}
\newcommand{\ctg}{{\rm cot}}
\newcommand{\degg}{{\mbox{deg}}}
 \newcommand{\Ad}{{\mbox{Ad}}}
\newcommand{\ad}{{\mbox{ad}}}
\newcommand{\codim}{{\mbox{codim}}}
\newcommand{\Gr}{{\mbox{Gr}}}
\newcommand{\coker}{{\rm coker}}
\newcommand{\id}{{\mbox{id}}}
\newcommand{\ord}{{\mbox{\rm ord}}}
\newcommand{\nat}{{\bf  N}}
\newcommand{\supp}{{\mbox{supp}}}
\newcommand{\spec}{{\mbox{spec}}}
\newcommand{\Ann}{{\mbox{Ann}}}
\newcommand{\aca}{{\aaaa_\C^\ast}}
\newcommand{\ck}{{\cal K}}
\newcommand{\tck}{{\tilde{\ck}}}
\newcommand{\tnk}{{\tilde{\ck}_0}}
\newcommand{\ceep}{{{\cal E}(E)^\prime}}
\newcommand{\ncE}{{{}^\naaa\cE}}
\newcommand{\cB}{{\cal B}}
\newcommand{\hc}{{{\cal HC}(\gaaa,K)}}
\newcommand{\vsl}{{V_{\sigma_\lambda}}}
\newcommand{\czg}{{\cZ(\gaaa)}}
\newcommand{\csl}{{\chi_{\sigma,\lambda}}}
\newcommand{\cR}{{\cal R}}
\def\hB{\hspace*{\fill}$\Box$ \newline\noindent}
\newcommand{\proof}{{\em Proof.$\:\:$}}
\newcommand{\varho}{\varrho}
\newcommand{\ind}{{\rm index}}
\newcommand{\Ind}{{\rm Ind}}
\maketitle
\mbox{}\\[0.5cm]
\newtheorem{prop}{Proposition}[section]
\newtheorem{lem}[prop]{Lemma}
\newtheorem{ddd}[prop]{Definition}
\newtheorem{theorem}[prop]{Theorem}
\newtheorem{kor}[prop]{Corollary}
\newtheorem{ass}[prop]{Assumption}
\newtheorem{con}[prop]{Conjecture}
\renewcommand{\imath}{{\it i}}

\tableofcontents
\mbox{}\\[0.5cm]
\section{Introduction}

Let $\Gamma\subset SL(2,\R)=:G$ be  a discrete subgroup
acting freely on the real hyperbolic plane $H^2=SL(2,\R)/S^1=:X$, i.e.,
the upper half space $\R^2_+$ with the metric $ds^2=\frac{1}{y^2}(dx^2+dy^2)$.
We assume this action to be geometrically finite and convex co-compact.
Under these circumstances $\Gamma$ will be called Fuchsian of the second kind.
Another model for the hyperbolic plane is the Poincar\'e disc $\{|z|<1\}\subset
\C$.
The corresponding group of isometries is $SU(1,1)$ which is, of course,
ismorphic to $SL(2,\R)$.

The group $G$ acts on the circle $S^1$
which can be identified with the
boundary $\partial X$ of $X$ using the Poincar\'e disc  model.
The set of accumulation points (in the Euclidean topology) of an arbitrary
orbit $\Gamma x$,
$x\in X$, is called the limit set $\Lambda\subset \partial X$.
The complement $\Omega:=\partial X\setminus \Lambda$ is called the ordinary
set.
$\Gamma$ acts freely and co-compactly on $\Omega$ with the compact quotient
$B$.
As a manifold $B$ is a finite union of circles.

Let $T\rightarrow S^1$ be the complexified tangent bundle of $S^1$.
It is $G$-homogeneous and we can form complex powers $T^\lambda\rightarrow
S^1$, $\lambda\in\C$.
Let $\partial $ be the fundamental vector field of the action of $U(1)\subset
G=SU(1,1)$.
It provides a natural trivialization of $T$.
Frequently we will identify a section $\phi
\partial^\lambda\in\Gamma(T^\lambda)$ with the function
$\phi$.
The  number $\lambda$ parametrizes a  principal series representation
$H^\lambda$ of $G$ on the Hilbert space
$L^2(S^1,T^{\lambda-1/2})$.
By $H^\lambda_{-\omega}$ we denote the space of its hyperfunction vectors.
As a topological vector space $H^\lambda_{-\omega}$
can be identified with the space of hyperfunction sections of
$T^{\lambda-1/2}$.
Thus $f\in H^\lambda_{-\omega}$ is a continuous functional on
$H^{-\lambda}_\omega:=C^\omega(S^1,T^{-\lambda-1/2})$.
\begin{ddd}
By $H^\lambda_{-\omega,\Lambda}$ we denote the subspace
of $H^\lambda_{-\omega}$ of hyperfunction
sections with support in the limit set $\Lambda$.
\end{ddd}
Since $\Lambda$ is a $\Gamma$-invariant closed subset
$H^\lambda_{-\omega,\Lambda}$
is a closed $\Gamma$-invariant subspace of $H^\lambda_{-\omega}$.

The goal of the present paper is to study the cohomology groups
$H^\ast(\Gamma,H^\lambda_{-\omega,\Lambda}) $.
Our interest in these cohomology groups is motivated by a conjecture
of S. Patterson \cite{patterson93} (see also \cite{bunkeolbrich947})
relating their dimensions with the order of the singularities
of the Selberg zeta function associated to $\Gamma$.
As a by-product we describe an easy way to obtain a meromorphic continuation
of the scattering matrix and the Eisenstein series. Here we recover  results of
Patterson \cite{patterson75}, \cite{patterson76}, \cite{patterson761}.

With the exception of  Section \ref{uio} the methods
extend to spherical principal series representations in
higher dimensions. To obtain the results of
Section \ref{uio} and to cover also non-spherical principal
series representations is still a challenging project.

{\bf Acknowledgment :}
{ \em We thank S. Patterson for explaining us his view on the topic
of the present paper. The results of Section \ref{uio} were obtained
while the authors visited the SFB 170 in G\"ottingen.}

\section{The case $\Ree(\lambda)>0$}

Let $A_\lambda:=\Delta-1/4+\lambda^2$, where
$\Delta$ is the Laplace-Beltrami operator of $X$.
Let $\cE=C^\infty(X)$ and $\cE(A_\lambda)=\ker(A_\lambda)\subset \cE$.
Helgason (\cite{helgason70}, \cite{helgason84} Introduction Thm. 4.3)  proved
that
for $\Ree(\lambda)>0$  the Poisson transform
$P_\lambda:H^\lambda_{-\omega}\rightarrow \cE(A_\lambda)$
is a $G$-isomorphism.
We characterize the subspace $P_\lambda(H^\lambda_{-\omega,\Lambda})\subset
\cE(A_\lambda)$
using certain semi-norms.
For all $W\subset X$  contained in finitely many translates of a fundamental
domain  of $\Gamma$ and
$k\in \nat_0$ we define
$$q_{W,k}(f):=\|A_\lambda^kf\|_{L^2(W)}\ .$$
By $\cE_\Lambda$ we denote the Frechet space of all functions $f\in\cE$ with
$q_{W,k}(f)<\infty$ for all such $W$ and $k\ge 0$.
Let $\cE_\Lambda(A_\lambda)=\cE_\Lambda\cap\ker(A_\lambda)$.
\begin{lem}\label{volley}
$$P_\lambda(H^\lambda_{-\omega,\Lambda}) =\cE_\Lambda(A_\lambda)\ .$$
\end{lem}
\proof
We first employ the upper half-plane model and assume that $\infty\in\Lambda$.
The Poisson transform $P_\lambda$ is given by the
kernel $P_\lambda((x,y),b):=(\frac{y}{y^2+(x-b)^2})^{\lambda+1/2}$,
where $b\in \R$.
The closures of the sets $W$ in the Euclidean topology are well separated from
$\Lambda$.
It is easy to see that for all $W$ the kernel
$b\to P_\lambda(.,b)$ defines an analytic function from a neighbourhood
of $\Lambda$ to $L^2(W)$. It also satisfies $A_\lambda P_\lambda(.,b)=0$.
If we pair the Poisson kernel with a hyperfunction supported on $\Lambda$ we
end up with
a smooth eigenfunction of $A_\lambda$ which is in $L^2(W)$ for all $W$. It
follows that
$$P_\lambda(H^\lambda_{-\omega,\Lambda}) \subset \cE_\Lambda(A_\lambda)\ .$$
It remains to show that the inverse (up to a scalar)  of $P_\lambda$,
the boundary value map $\beta_\lambda$, satisfies
$$\beta_\lambda(\cE_{\Lambda}(A_\lambda))\subset H^\lambda_{-\omega,\Lambda}\
.$$
For this argument it is more convenient to employ the ball model.
We use polar coordinates $(r,\alpha)$, $r\in [0,1)$, $\alpha\in (0,2\pi]$
 of $X$. Let $f\in \cE_{\Lambda}(A_\lambda)$ and $\phi\in H_\omega^{-\lambda}$.
Then
$$\langle \beta_\lambda(f),\phi\rangle = \lim_{r\to 1} (1-r^2)^{\lambda-1/2}
\int_{S^1} f(r\ee^{\imath\alpha}) \phi(\ee^{\imath\alpha})
\frac{d\alpha}{2\pi}\ .$$
Now let $(a,b)\subset \Omega$. The $L^2$-condition satisfied by $f$ along
$\Omega$ allows us to restrict the integration to the complement of $(a,b)$,
i.e.,
$$\langle \beta_\lambda (f),\phi\rangle = \lim_{r\to 1} (1-r^2)^{\lambda-1/2}
\int_{S^1\setminus (a,b)} f(r\ee^{\imath\alpha}) \phi(\ee^{\imath\alpha})
\frac{d\alpha}{2\pi}\ .$$
Thus $\beta_\lambda(f)$ is in fact a continuous functional on germs of analytic
sections of
$T^{-\lambda-1/2}$ on $S^1\setminus (a,b)$, i.e,  $\beta_\lambda(f)\in{\cal
A}^\prime(S^1\setminus (a,b),T^{ \lambda-1/2})$. The hyperfunction sections of
$T^{\lambda-1/2}$ on
$(a,b)$ are ${\cal B}((a,b),T^{ \lambda-1/2})\cong {\cal A}^\prime(S^1,T^{
\lambda-1/2})/{\cal A}^\prime(S^1\setminus (a,b),T^{ \lambda-1/2})$.
Hence the restriction to $(a,b)$ of the  hyperfunction $ \beta_\lambda(f)$
vanishes.
Since this holds for all intervals $(a,b)\subset \Omega$ we conclude
$\supp \:\beta_\lambda(f)\subset \Lambda$.
The lemma follows. \hB

The proof of the following lemma is completely analogous to that of
\cite{bunkeolbrich947},
Lemma 2.4.
\begin{lem}\label{cri}
  $\cE_\Lambda$ is $\Gamma$-acyclic.
\end{lem}
The following lemma shows that
\begin{equation}\label{nummer1}
0\longrightarrow
H^\lambda_{-\omega,\Lambda}\stackrel{P_\lambda}{\longrightarrow}\cE_\Lambda
\stackrel{A_\lambda}{\longrightarrow}\cE_\Lambda\longrightarrow 0
\end{equation}
is a $\Gamma$-acyclic resolution of $H^\lambda_{-\omega,\Lambda}$.
\begin{lem}\label{wolf}
$A_\lambda:\cE_\Lambda\rightarrow \cE_\Lambda$ is surjective.
\end{lem}
\proof
Let ${}^tA_\lambda:\cE_\Lambda^\prime\rightarrow \cE_\Lambda^\prime$
be the dual operator. It is enough to show that ${}^tA_\lambda$
is injective and has  closed range (see Treves \cite{treves67}
for basic techniques to deal with  surjectivity questions).

Since $C_c^\infty(X)$ is dense in $\cE_\Lambda$ we can embed
$\cE_\Lambda^\prime$
into the distributions on $X$. If $f\in \cE_\Lambda^\prime$, then there
exist  $k,W$ such that  $|\langle f,\phi\rangle|\le C q_{W,k}(\phi)$ for all
$\phi\in C_c^\infty(X)$. Hence $f$ is a distribution of order at most $2k$ with
support on $W$. The bound of the order follows from the fact that
$q_{W,k}(\phi)$
estimates the local $H^{2k }$-norm of $\phi$.

If $f\in \cE_\Lambda^\prime$, then it is supported on a finite number of
fundamental
domains. If in addition ${}^tA_\lambda f=0$, then $f=0$, since $f$  is real
analytic
and vanishes on an open set. This shows the injectivity of ${}^tA_\lambda$.

The range of ${}^tA_\lambda$ is closed, if $\im({}^tA_\lambda)\cap B^\prime$
is closed for any bounded set $B^\prime\subset \cE_\Lambda^\prime$.
Let $B^\prime$ be bounded. Since $\cE_\Lambda$ is Fr\'echet, by a theorem of
Schwartz
\cite{schwartz66}, p.90, there are   $W$, $k$, $\epsilon>0$
such that $X\setminus W$ is connected and  $|\langle f,\phi\rangle| \le 1$ for
all $f\in B^\prime$, $\phi\in C_c^\infty(X)$ with
$q_{W,k}(\phi)\le \epsilon$. It follows that  $\supp(f)\subset W$, $\forall
f\in B^\prime$.

Let ${}^tA_\lambda h_i=f_i\in B^\prime$ converge to $f$.
Unique continuation for ${}^tA_\lambda$  implies that $\supp (h_i)\subset W$.
We extend
$\Delta-1/4$  to a non-negative self-adjoint unbounded operator
on $L^2(X)$.
Since either $\Imm(\lambda^2)\not=0$
or $\lambda^2>0$ the operator ${}^tA_\lambda=\Delta-1/4+ \lambda^2 $
has a bounded inverse on $L^2(X)$.
We use ${}^tA_\lambda$ in order to define the scale of Sobolev spaces
$H^l(X)$, $l\in 2\Z$. Here $A_\lambda:H^l(X)\rightarrow H^{l-2}(X)$
is an isometry and $H^0(X)=L^2(X)$.
It is easy to see that $f_i$ converges in the direct limit $\lim_{l\to\infty
}H^{-l}(X)$.
In fact it already converges in $H^{-2k}(X)$.
Hence $h_i\to h={}^tA_\lambda^{-1}(f)$. Since $\supp(h_i)
\subset W$ we also have $\supp(h)\subset W$. It follows that $h\in B^\prime$
and ${}^tA_\lambda h=f$.
This shows the closedness of the range of ${}^tA_\lambda$.
The lemma follows. \hB

In order to calculate the cohomology  of $\Gamma$ with coefficients in
$H^\lambda_{-\omega,\Lambda}$
we take the subcomplex of $\Gamma$-invariants of the complex (\ref{nummer1}).
Let $Y=\Gamma\backslash X$ be the Riemann surface corresponding to $\Gamma$.
Let $A_{\lambda,Y}=\Delta_Y-1/4+\lambda^2$, where $\Delta_Y$ is the Laplacian
of $Y$.
Let $\cE_Y:=\{f\in C^\infty(Y)\ |\ A_{\lambda,Y}^kf\in L^2(Y)\ \forall
k\in\nat_0 \}$.
The cohomology $H^\ast(\Gamma,H^\lambda_{-\omega,\Lambda})$
is the cohomology of
$$0\longrightarrow
\cE_Y\stackrel{A_{\lambda,Y}}{\longrightarrow}\cE_Y\longrightarrow 0.$$
Since $-\lambda^2$ is separated from the essential spectrum $[0,\infty)$   of
$\Delta_Y-1/4$
(see \cite{patterson75}) we obtain
\begin{prop}\label{erst}
If $\Ree(\lambda)>0$,
then $H^p(\Gamma,H^\lambda_{-\omega,\Lambda})=0$ for all $p\ge 2$  and
$H^0(\Gamma,H^\lambda_{-\omega,\Lambda})=
H^1(\Gamma,H^\lambda_{-\omega,\Lambda})= \ker_{L^2}(A_{\lambda,Y})$. There is a
finite set
$\Sigma\subset (0,\infty)$ such that $\ker_{L^2}(A_{\lambda,Y})$   is
non-trivial iff $\lambda\in  \Sigma\subset \C$. Here $ \ker_{L^2}(A_\lambda)$
is always finite-dimensional.
\end{prop}
The finiteness of the point spectrum of $\Delta_Y$ was shown e.g. in
Lax/Phillips \cite{laxphillips82}.

\section{Intertwining operators and the scattering matrix}\label{sss}

We introduce the spaces $H^\lambda_\infty(\Omega)$ and
$H^\lambda_{-\omega,(\Lambda)}$ of smooth functions on $\Omega$
and hyperfunctions on $S^1$ with singular support in $\Lambda$, respectively.
The index $\lambda$ indicates the way $\Gamma$ acts on these spaces.
Since the sheaf of hyperfunctions is flabby we have an exact sequence
$$0\longrightarrow H^\lambda_{-\omega,\Lambda}\longrightarrow
H^\lambda_{-\omega,(\Lambda)}\stackrel{res}{\longrightarrow}
H^\lambda_\infty(\Omega)\longrightarrow 0$$
of $\Gamma$-modules.
The long exact cohomology sequence gives
\begin{eqnarray*}
&&0\longrightarrow {}^\Gamma H^\lambda_{-\omega,\Lambda}\longrightarrow
{}^\Gamma H^\lambda_{-\omega,(\Lambda)}\stackrel{res}{\longrightarrow}
{}^\Gamma H^\lambda_\infty(\Omega)\longrightarrow\\
&&\longrightarrow H^1(\Gamma,H^\lambda_{-\omega,\Lambda})\longrightarrow\dots\
{}.
\end{eqnarray*}
Recall the definition of the finite set
$\Sigma:=\{\lambda\in\C\ | 1/2 > \Ree(\lambda)>0,
H^0(\Gamma,H^\lambda_{-\omega,\Lambda})\not=0 \}$.
For all $\lambda\in\C\setminus\Sigma $, $\Ree(\lambda)>0$, by Proposition
\ref{erst}
an invariant  function $\phi_\lambda\in {}^\Gamma H^\lambda_\infty(\Omega)$
has a unique invariant extension as a hyperfunction $\Phi_\lambda\in {}^\Gamma
H^\lambda_{-\omega,(\Lambda)}$, i.e.,  $ res(\Phi_\lambda)=\phi_\lambda$.

We fix a smooth positive invariant function $\varho\in {}^\Gamma
H^{3/2}_\infty(\Omega)$
(for existence see Perry \cite{perry89}).
Multiplication by
$$\varho^s:H^\lambda_\infty(\Omega)\rightarrow H^{\lambda+s}_\infty(\Omega)$$
is an isomorphism of $\Gamma$-modules.
Let $\phi\in {}^\Gamma H^{1/2}_\infty(\Omega)$ and set
$\phi_\lambda:=\varho^{\lambda-1/2}\phi\in {}^\Gamma H^\lambda_\infty(\Omega)$.
Then $\phi_\lambda$ is a holomorphic family of smooth functions
on $\Omega$. For $\lambda\in \C\setminus \Sigma$, $\Ree(\lambda)>0$,
let $\Phi_\lambda\in {}^\Gamma H^\lambda_{-\omega,(\Lambda)}$ be
the unique extension of $\phi_\lambda$.
\begin{lem}\label{aass}
For $\lambda\in \C\setminus \Sigma$, $\Ree(\lambda)>0$, the family
$\Phi_\lambda$ is  a holomorphic family of hyperfunctions.
\end{lem}
\proof
Consider the Poincar\'e series
$$\sum_{g\in\Gamma}  g^\prime(b)^{\lambda+1/2}\ $$
for some $b\in\Omega$.
Here $g^\prime(b)$ is the conformal dilatation of the map $b\mapsto g(b)$ with
respect
to the Euclidean metric on $S^1=\partial X$ in the Poincar\'e disc model.
Let $\delta_\Gamma \in\R$ be the smallest number such that the Poincar\'e
series
converges  for all $\lambda\in\C$ with $\Ree(\lambda)>\delta_\Gamma$.
For $ \Ree(\lambda)>\delta_\Gamma$  the convergence is independent of the
choice of $b\in\Omega$ and uniform
in compact subsets of $\Omega$.
It is known that $\delta_\Gamma < 1/2$ and that $\delta_\Gamma+1/2$
is the Hausdorff dimension of $\Lambda$ (Patterson \cite{patterson762},
Sullivan \cite{sullivan79}).
We claim that for $\Ree(\lambda)>\delta_\Gamma$
the hyperfunction extension $\Phi_\lambda$ of $\phi_\lambda$ is  given by
integration against $\phi_\lambda$.
This extension is in fact a measure.
Let $f\in C(S^1)$ and  $F\subset\Omega$ be a compact fundamental domain for
$\Gamma$. We have
\begin{eqnarray}
<\Phi_\lambda,f\partial^{-\lambda-1/2}>&:=&\int_\Omega
\varho^{\lambda-1/2}(b)\phi(b) f(b) db\label{ppoo}\\
&=&\sum_{g\in\Gamma} \int_{g^{-1}F} \varho^{\lambda-1/2}(b)\phi(b) f(b)
db\nonumber\\
&=&\sum_{g\in\Gamma} \int_{F} \varho^{\lambda-1/2}(gb)\phi(gb) f(gb)
d(gb)\nonumber\\
&=&\sum_{g\in\Gamma} \int_{F} g^\prime(b)^{\lambda+1/2}
\varho^{\lambda-1/2}(b)\phi(b) f(gb) db\nonumber\\
&=&\int_{F}(\sum_{g\in\Gamma}g^\prime(b)^{\lambda+1/2} f(gb))
\varho^{\lambda-1/2}(b)\phi(b)db
\nonumber\ .
\end{eqnarray}
Since $f$ is uniformly bounded on $S^1$ the right-hand side converges for
$\Ree(\lambda)>\delta_\Gamma$.
It is easy to see that   $\Phi_\lambda$ is  a holomorphic family of measures
with respect to the parameter $\lambda$.  Hence it is a holomorphic family of
hyperfunctions for $\Ree(\lambda)>\delta_\Gamma$.

It remains to show that $\Phi_\lambda$ is holomorphic in the strip
$0<\Ree(\lambda)<1/2$,
$\lambda\not\in\Sigma$.
This requires analysis on the symmetric space $X$ and its quotient
$Y=\Gamma\backslash X$.
Let $T:(0,\epsilon)\times \Omega\rightarrow X$, $\epsilon>0$,
be the coordinates of a $\Gamma$- invariant neighbourhood of
$\Omega$ introduced by Perry \cite{perry87}.
These coordinates depend on the choice of $\varho$
and are defined as follows:
\begin{eqnarray*}
y(r,b)&=&\frac{r\varho(b)}{1+\frac{r^2}{4}\varho^\prime(b)^2}\\
x(r,b)&=&b-\frac{1}{2}\frac{r^2
\varho(b)\varho^\prime(b)}{1+\frac{r^2}{4}\varho^\prime(b)^2}\ .
\end{eqnarray*}
Here $(x,y)$, $y>0$, are coordinates of the upper half-plane model.
The map $T$ satisfies $gT(r,b)=T(r,gb)$ for $g\in\Gamma$, $r\in(0,\epsilon)$,
$b\in\Omega$.
Hence we obtain induced coordinates of a collar neighbourhood $E\subset Y$ at
infinity:

\begin{equation}\label{coor}T_Y:(0,\epsilon)\times B \rightarrow Y\
.\end{equation}
In these coordinates the Laplacian $\Delta_Y$ has the form
$$\Delta_Y=-(r\frac{\partial}{\partial r})^2+r\frac{\partial}{\partial r}+r^2
\Delta_B +r P(r\frac{\partial}{\partial r},r \frac{\partial}{\partial b})\ ,$$
where $\Delta_B$ is the Laplacian on $B$ with respect to the metric induced by
$\varho^{-2}(db)^2$
and $P$ is of second order with bounded coefficients.
Let $A_{\lambda,Y}=\Delta_Y-1/4+\lambda^2$.
Asymptotically the volume form is $r^{-2}$-times a bounded form.
Fix a cut-off function $\chi(r)$ being zero near $r=\epsilon$ and one near
$r=0$.
Form $f_\lambda(r,b)=\chi(r)r^{-\lambda+1/2}\phi(b)$.
One can check that $g_\lambda:=A_{\lambda,Y}f_\lambda$
is a holomorphic family in $L^2(Y)$ for $\Ree(\lambda)<1$.
For $\Ree(\lambda)>0$, $\lambda\not\in\Sigma$,  we can solve
$$A_{\lambda,Y}h_\lambda=-g_\lambda$$
in $L^2(Y)$
and obtain a holomorphic family of eigenfunctions
$F_\lambda:=f_\lambda+h_\lambda\in \ker(A_{\lambda,Y})$.
Lifting $F_\lambda$ to $X$ we find a holomorphic family of invariant
eigenfunctions
$\tilde{F}_\lambda$.
 The boundary value $\beta_\lambda(\tilde{F}_\lambda)$ is a hyperfunction and
the desired extension $\Phi_\lambda$ of $\phi_\lambda$ for $\Ree(\lambda)>0$,
$\lambda\not\in\Sigma$.
In fact on sets $W$ we have $\tilde{F}_\lambda(x,y)\sim
y^{-\lambda+1/2}\varho^{\lambda-1/2}(x)\phi(x)$
(mod $L^2(W)$) and thus $res\:\beta_\lambda(\tilde{F}_\lambda)=\phi_\lambda$.
For $\delta_\Gamma<\Ree(\lambda)<1$ the extension just constructed coincides
with
the extension defined by (\ref{ppoo}) by uniqueness.
Since the boundary value map depends holomorphically on $\lambda$ for
$\Ree(\lambda)>0$   the lemma follows.
\hB
 \begin{ddd}
We define the Eisenstein series associated to $\phi\in {}^\Gamma
H^{1/2}_\infty(\Omega)$ and $\lambda\in \C$ with   $\Ree(\lambda)>0$,
$\lambda\not\in\Sigma$, by
$$E(\phi,\lambda)=P_\lambda(\Phi_\lambda)\ .$$
\end{ddd}
The term 'series' stems from the fact that for $\Ree(\lambda)>\delta_\Gamma$
it can be written as
$$E(\phi,\lambda)(p)=\sum_{g\in\Gamma} \int_F P_\lambda(gp,b)\phi_\lambda(b)
db$$
(compare \cite{patterson75}, \cite{mandouvalos88}, \cite{perry89},
\cite{patterson89}),
since the Lebesgue measure of the limit set is zero and we can alternatively
integrate over
all of $S^1$:
$$E(\phi,\lambda)(p)= \int_{S^1 }P_\lambda(p,b)\Phi_\lambda(b) db\ .$$

We now discuss the singularities of the extension $\Phi_\lambda$ at
$\lambda_0\in\Sigma$. Let  $\{F_i\}$, $i=1,\dots,k$,
be an orthonormal base of  $ \ker_{L^2}(A_{\lambda_0,Y})$. Then near
$\lambda_0$
the family $F_\lambda$ has the form
$$F_\lambda=\frac{1}{\lambda_0^2-\lambda^2}\sum_{i=1}^k \langle F_i,g_\lambda
\rangle_{L^2(Y)} F_i+\hat{F}_\lambda\ ,$$
where $\hat{F}_\lambda$ is holomorphic and orthogonal to $
\ker_{L^2}(A_{\lambda_0,Y})$.
Applying the boundary value map $\beta_\lambda$ we obtain
\begin{equation}\label{eees}\Phi_\lambda=\frac{1}{\lambda_0^2-\lambda^2}
\sum_{i=1}^k \langle F_i,g_\lambda \rangle_{L^2(Y)}
\beta_\lambda(\tilde{F_i})+\hat{\Phi}_\lambda\ ,\end{equation}
where $\hat{\Phi}_\lambda$ is holomorphic.
Note that $F_i(r,b)\sim r^{\lambda+1/2}\beta_Y(F_i)(b)+ o(r^{\lambda+1/2})$,
where $\beta_Y$ is the boundary value on $Y$.
In more detail, $\beta_Y(F_i)(.)=\lim_{r\to 0} r^{-\lambda-1/2}F_i(r,.)$, where
the limit exists in the $L^2(B)$-sense.
Applying Green's formula we obtain
\begin{eqnarray*}
\langle F_i,g_{\lambda_0} \rangle_{L^2(Y)}&=&\langle F_i,A_{{\lambda_0},Y}
f_{\lambda_0} \rangle_{L^2(Y)} - \langle A_{{\lambda_0},Y} F_i, f_{\lambda_0}
\rangle_{L^2(Y)} \\
&=&\lim_{r\to 0}\langle F_i,r\frac{d}{dr} f_{\lambda_0} \rangle_{L^2(\{r
\}\times B)}-\lim_{r\to 0}\langle r\frac{d}{dr}F_i,f_{\lambda_0}
\rangle_{L^2(\{r \}\times B)}   \\
&=&-2{\lambda_0} \langle \beta_Y(F_i) ,\phi\rangle_{L^2(B)}\ .
\end{eqnarray*}
Taking the Poisson transform of (\ref{eees}) and using
$$P_\lambda\circ\beta_\lambda=\frac{1}{\sqrt{2}}\:
\frac{\Gamma(\lambda)}{\Gamma(1/2+\lambda)}
\id$$ we obtain
\begin{prop}
The Eisenstein series $E(\phi,\lambda)$ is meromorphic on $\Ree(\lambda)>0$
with singularities in $\Sigma$ and residues
$$\frac{1}{\sqrt{2}}\:\frac{\Gamma(\lambda_0)}{
\Gamma(1/2+\lambda_0)}\sum_{i=1}^k  \langle \beta_Y(F_i) ,\phi\rangle_{L^2(B)}
F_i $$
at $\lambda=\lambda_0\in\Sigma$.
\end{prop}
Recall the  Knapp-Stein intertwining operator
$J_\lambda:H^\lambda_{-\omega}\rightarrow H^{-\lambda}_{-\omega}$.
For $\Ree(\lambda)< 0$ it is given by the integral operator
$$(J_\lambda\phi)(b)=\int_\R |b-b^\prime|^{-1-2\lambda}\phi(b^\prime)
db^\prime,\quad \phi\in C_c^\infty(R)\ .$$
For $\Ree(\lambda)\ge 0$ it is defined by analytic continuation.
$J_\lambda$ is a meromorphic family of elliptic pseudodifferential operators
with singularities at $\lambda\in\nat_0$.
The restriction of $J_\lambda$ to $H^\lambda_{-\omega,(\Lambda)}$
defines a meromorphic family of $\Gamma$-intertwining operators
$H^\lambda_{-\omega,(\Lambda)}\rightarrow H^{-\lambda}_{-\omega,(\Lambda)}$.
\begin{ddd}\label{scat}
For $\Ree(\lambda)>0$, $\lambda\not\in\Sigma$, the scattering matrix
$S_\lambda$ is the
operator  $$S_\lambda:{}^\Gamma H^\lambda_\infty(\Omega)\rightarrow {}^\Gamma
H^{-\lambda}_\infty(\Omega)$$
given by
$${}^\Gamma H^\lambda_\infty \ni\phi \mapsto \Phi\in {}^\Gamma
H^\lambda_{-\omega,(\Lambda)}\mapsto res\circ J_\lambda(\Phi)\in {}^\Gamma
H^{-\lambda}_\infty(\Omega)\ .$$
\end{ddd}
For $\Ree(\lambda)>\delta_\Gamma$, $\lambda\not\in\nat$, and $\phi\in {}^\Gamma
C^\infty(\Omega)$ we have
\begin{eqnarray*}
S_\lambda\phi_\lambda(b)&=&\int_{\Omega}
|b-b^\prime|^{-1-2\lambda}\varho^{\lambda-1/2}\phi(b^\prime) db^\prime\\
&=&\sum_{g\in\Gamma}\int_{g^{-1}
F}|b-b^\prime|^{-1-2\lambda}\varho^{\lambda-1/2}(b^\prime)\phi(b^\prime)
db^\prime\\
&=&\sum_{g\in\Gamma}\int_{F}|b-gb^\prime
|^{-1-2\lambda}\varho^{\lambda-1/2}(gb^\prime)\phi(gb^\prime) d(gb^\prime)\\
&=&\sum_{g\in\Gamma}\int_{F}
|b-gb^\prime|^{-1-2\lambda}g^\prime(b^\prime)^{\lambda+1/2}
\varho^{\lambda-1/2}(b^\prime)\phi(b^\prime) db^\prime\ ,
\end{eqnarray*}
where the singularity of the kernel near  $b=b^\prime$ has to be interpreted in
the
regularized sense.
Thus for $\Ree(\lambda)>\delta_\Gamma$, $\lambda\not\in\nat$, the scattering
matrix is given by the integral kernel
$$ \sum_{g\in\Gamma}
|b-gb^\prime|^{-1-2\lambda}g^\prime(b^\prime)^{\lambda+1/2}$$
and this coincides with other definitions in the literature
(\cite{patterson75}, \cite{patterson89}, \cite{perry89}).
The scattering matrix
is a meromorphic family of elliptic pseudodifferential operators on $B$ in the
following sense.
Let $\phi\in C^\infty(B)={}^\Gamma H^{1/2}_\infty(\Omega)$.
Define $S(\lambda):C^\infty(B)\rightarrow C^\infty(B)$ by
$$S(\lambda)\phi:=\varho^{1/2+\lambda}S_\lambda\varho^{\lambda-1/2}\phi\ .$$
\begin{prop}\label{form}
The family
$$S(\lambda)\phi:=\varho^{1/2+\lambda}S_\lambda\varho^{\lambda-1/2}\phi$$ is
meromorphic for $\Ree(\lambda)>0$ with poles
at $\Sigma$ and   residues
\begin{equation}\label{pppo}\phi\mapsto
\frac{1}{\sqrt{2}}\:\frac{\Gamma(\lambda_0)}{\Gamma(1/2+\lambda_0)}
\sum_{i=1}\langle \beta_Y(F_i),\phi\rangle_{L^2(B)}  \beta_Y(
F_i)\end{equation}
at $\lambda=\lambda_0\in\Sigma$.
Moreover, it has poles at $\lambda\in\nat$ (induced by the singularities of the
on-diagonal regularization of the intertwining operator)
and $S(k+1/2)=\varho^{k+1} d^{2k+1} \varho^{k}$, $k\in\nat_0$.
\end{prop}
Equation (\ref{pppo}) is obtained from (\ref{eees}) and
$$J_\lambda\circ\beta_\lambda=\frac{1}{\sqrt{2}}\:\frac{\Gamma(\lambda)}{
\Gamma(1/2+\lambda)}\beta_{-\lambda}\ .$$
In Section \ref{scats} we obtain a meromorphic continuation of the scattering
matrix and the Eisenstein series to all of $\C$.
\section{The case $\Ree(\lambda)<0$, $\lambda\not=-1/2,-3/2,\dots$}\label{s4}

For $2\lambda\not\in\Z$ the intertwining operator $J_\lambda$
is an isomorphism. By restriction we obtain an isomorphism of
$\Gamma$-modules $J_\lambda:H^\lambda_{-\omega,(\Lambda)}\rightarrow
H^{-\lambda}_{-\omega,(\Lambda)}$. The operator $J_{ \lambda}$ has first order
poles at $\lambda=k\in  \nat$.  At these points we consider its renormalization
$\tilde{J}_{k}:=\lim_{\lambda\to k} (\lambda-k)J_{\lambda}$.
It satisfies
$\tilde{J}_{k}\circ J_{-k}=-\frac{\pi}{2k}\id$ and provides an isomorphism
$\tilde{J}_{k}:H^k_{-\omega,(\Lambda)}\rightarrow H^{-k}_{-\omega,(\Lambda)}$ .

\begin{lem}\label{asd}
 $H^\lambda_{-\omega,(\Lambda)}$ is $\Gamma$-acyclic.
\end{lem}
\proof
It is enough to show that $H^{-\lambda}_{-\omega,(\Lambda)}$
is $\Gamma$-acyclic.
Since $\Gamma$ acts freely and properly on $\Omega$
the modules $H^{\lambda}_\infty(\Omega)$ are $\Gamma$-acyclic
by  \cite{bunkeolbrich947}, Lemma 2.4, for all $\lambda\in\C$.
Consider the exact sequence of $\Gamma$-modules
$$0\longrightarrow H^{-\lambda}_{-\omega,\Lambda}\longrightarrow
H^{-\lambda}_{-\omega,(\Lambda)}\stackrel{res}{\longrightarrow}
H^{-\lambda}_\infty(\Omega)\longrightarrow 0\ .$$
The associated long exact sequence is
\begin{eqnarray*}
\lefteqn{0\longrightarrow
H^0(\Gamma,H^{-\lambda}_{-\omega,\Lambda})\longrightarrow
H^0(\Gamma,H^{-\lambda}_{-\omega,(\Lambda)})\stackrel{res}{\longrightarrow}
H^0(\Gamma,H^{-\lambda}_\infty(\Omega))\stackrel{\delta}{
\longrightarrow}}\hspace{2.5cm}\\
&& \longrightarrow H^1(\Gamma,
H^{-\lambda}_{-\omega,\Lambda})\longrightarrow
H^1(\Gamma,H^{-\lambda}_{-\omega,(\Lambda)}) \longrightarrow 0\ .
\end{eqnarray*}
In case that $-\lambda\not\in\Sigma$ we immediately deduce
$H^1(\Gamma,H^{-\lambda}_{-\omega,(\Lambda)})=0$.
Thus assume that $-\lambda\in\Sigma$. We have to show that $\delta$
is surjective. Let $d:=\dim {}^\Gamma H^{-\lambda}_{-\omega,\Lambda} =\dim
H^1(\Gamma,H^{-\lambda}_{-\omega,\Lambda})$. It is suffices to show that
$$\dim\coker(res :{}^\Gamma H^{-\lambda}_{-\omega,(\Lambda)} \rightarrow
{}^\Gamma H^{-\lambda}_\infty(\Omega)) \ge d\ .$$
Let $\tilde{W}\subset{}^\Gamma H^{\lambda}_\infty(\Omega)$ be the subspace
spanned
by the boundary values $\varrho^{\lambda-1/2}\beta_Y(F)$ for all $F\in
\ker_{L^2}(A_{-\lambda,Y})$
and $V\subset {}^\Gamma H^{-\lambda}_\infty(\Omega)$ be the orthogonal
complement of $\tilde{W}$.  Then $\codim\:V = d$. It is sufficient to show that
$res ({}^\Gamma H^{-\lambda}_{-\omega,(\Lambda)})  \subset V$.
Let $\Phi\in {}^\Gamma H^{-\lambda}_{-\omega,(\Lambda)}$.
We must verify that $\langle \beta_Y(F) , \varrho^{ \lambda +1/2}
res(\Phi)\rangle_{L^2(B)}=0$ for all $F\in  \ker_{L^2}(A_{-\lambda,Y})$.
Thus let $E=P_{-\lambda}\Phi\in C^{\infty}(Y)$.
Then by Green's formula
\begin{eqnarray*}
\lefteqn{-2{\lambda } \langle \beta_Y(F) , \varrho^{ \lambda +1/2}res(\Phi)
\rangle_{L^2(B)}}\\
&=&\lim_{r\to 0}\langle F,r\frac{d}{dr} E \rangle_{L^2(\{r \}\times
B)}-\lim_{r\to 0}\langle r\frac{d}{dr}F,E \rangle_{L^2(\{r \}\times B)}   \\
&=&\langle F,A_{{\lambda },Y} E \rangle_{L^2(Y)} - \langle A_{{\lambda },Y} F ,
E \rangle_{L^2(Y)} \\
&=&0\ .
\end{eqnarray*}
This finishes the proof of the Lemma for $-\lambda\in\Sigma$.
 \hB
If $\lambda\not\in -\nat_0$, then $J_{-\lambda}\circ J_\lambda=-\frac{\cot(\pi
\lambda)}{2\lambda}\id$.
By Lemma \ref{asd} we obtain the following $\Gamma$-acyclic resolution of
$H^\lambda_{-\omega,\Lambda}$:
$$0\longrightarrow
H^\lambda_{-\omega,\Lambda}\stackrel{
J_\lambda}{\longrightarrow}H^{-\lambda}_{-\omega,
(\Lambda)}
\stackrel{res\circ J_{-\lambda}}{\longrightarrow}
H^{\lambda}_\infty(\Omega)\longrightarrow 0\ .$$
Restricting to the $\Gamma$-invariants and taking
Definition \ref{scat} of the scattering matrix  into account we obtain
\begin{prop}
For $\Ree(\lambda)<0$, $\lambda\not\in-\frac{1}{2}\nat_0$,
$-\lambda\not\in\Sigma$,
we have
\begin{eqnarray*}
H^p(\Gamma,H^\lambda_{-\omega,\Lambda})&=&0, \quad p\ge 2\ ,\\
H^0(\Gamma,H^\lambda_{-\omega,\Lambda})&=&\ker(S_{-\lambda})\ ,\\
H^1(\Gamma,H^\lambda_{-\omega,\Lambda})&=&\coker(S_{-\lambda})\ .
\end{eqnarray*}
Since $S_{-\lambda}$ is elliptic and $\ind(S_{-\lambda})=0$  we have
$$\dim H^0(\Gamma,H^\lambda_{-\omega,\Lambda})  = \dim
H^1(\Gamma,H^\lambda_{-\omega,\Lambda})<\infty\ .$$
\end{prop}
Now we discuss the points $\lambda$, where the intertwining operator  or
the scattering matrix have singularities.
First we consider the points $k\in\nat$.
We obtain the $\Gamma$-acyclic resolution
$$0\longrightarrow H^k_{-\omega,\Lambda}\stackrel{J_{-k}}{
\longrightarrow}H^{-k}_{-\omega,(\Lambda)}\stackrel{res\circ
\tilde{J}_{k}}{\longrightarrow} H^{k}_\infty(\Omega)\longrightarrow 0\ .$$
 Let $\tilde{S}_{k}:=\lim_{\lambda\to k}(\lambda-k)S_{\lambda}$ be the
renormalized scattering matrix defined in
a similar way as $\tilde{J}_k$.
Then we have
\begin{prop}
For $k\in\nat$, $k\not=0$, we have
\begin{eqnarray*}
H^p(\Gamma,H^{-k}_{-\omega,\Lambda})&=&0, \quad p\ge 2\ ,\\
H^0(\Gamma,H^{-k}_{-\omega,\Lambda})&=&\ker(\tilde{S}_{k})\ ,\\
H^1(\Gamma,H^{-k}_{-\omega,\Lambda})&=&\coker(\tilde{S}_{k})\ .
\end{eqnarray*}
Since $\tilde{S}_{-\lambda}$ is elliptic and  $\ind(\tilde{S}_{k})=0$  we have
$$\dim H^0(\Gamma,H^{-k}_{-\omega,\Lambda})  = \dim
H^1(\Gamma,H^{-k}_{-\omega,\Lambda})<\infty\ .$$
\end{prop}
We now consider the case $-\lambda\in\Sigma$.
By Lemma \ref{asd} we have $\dim\coker(res :{}^\Gamma
H^{-\lambda}_{-\omega,(\Lambda)} \rightarrow {}^\Gamma
H^{-\lambda}_\infty(\Omega)) = d$.
In fact $res ({}^\Gamma H^{-\lambda}_{-\omega,(\Lambda)}) =V$ is the orthogonal
complement
of the space spanned by   the boundary values of the elements of $
\ker_{L^2}(A_{-\lambda,Y})$ and thus $V$ is closed.
We choose a complement $W\subset {}^\Gamma H^{-\lambda}_\infty(\Omega)$ to $V$.
By Proposition \ref{form} the singular part of the scattering matrix
$S_{-\lambda}$
vanishes on $V$.
We fix an identification $i:W\stackrel{\sim}{\rightarrow}
H^{-\lambda}_{-\omega,\Lambda} $
 and   define a regularized scattering matrix by
$(S^{reg}_{-\lambda})_{|V}=S_{-\lambda}$ and
$(S^{reg}_{-\lambda})_{|W}=res\circ J_{-\lambda}\circ i$. Then
$S^{reg}_{-\lambda}$ is still an elliptic pseudodifferential
operator of index zero.
For any $\phi\in V$ by (\ref{eees}) we can construct the extension $\Phi\in
{}^\Gamma H^{-\lambda}_{-\omega,(\Lambda)}$
thus obtaining a map $ext:V\rightarrow {}^\Gamma
H^{-\lambda}_{-\omega,(\Lambda)}$
which is right inverse to $res$. It induces a decomposition
${}^\Gamma H^{-\lambda}_{-\omega,(\Lambda)}=
\im(ext)\oplus {}^\Gamma H^{-\lambda}_{-\omega,\Lambda}$.
Then
$$\tilde{res}= res_{|\im(ext)}\oplus i^{-1}  :\im(ext)\oplus
H^{-\lambda}_{-\omega,\Lambda}\stackrel{\sim}{\rightarrow} V\oplus W$$
satisfies $S^{reg}_{-\lambda}\circ\tilde{res}=res\circ J_{-\lambda} $.
Hence we can identify
$H^i(\Gamma,H^{\lambda}_{-\omega,\Lambda})$ for $i=0,1$ with the kernel and
cokernel
of $S^{reg}_{-\lambda}$, respectively.
We claim that $S^{reg}_{-\lambda}(V)$ is transverse to $res\circ
J_{-\lambda}({}^\Gamma H^{-\lambda}_{-\omega,\Lambda})$.
We employ the meromorphic continuation of the scattering matrix to all of $\C$
and
its functional equation obtained in Section \ref{scats}. Using Lemma
\ref{sctti} we find $res\circ J_{-\lambda}({}^\Gamma
H^{-\lambda}_{-\omega,\Lambda})\subset \ker(S_{\lambda})$ and $S_\lambda\circ
(S^{reg}_{-\lambda})_{|V}=-\frac{\cot(\pi \lambda)}{2\lambda} \id_V$.
Also, $res\circ J_\lambda\circ i^{-1}$ is injective. In fact a nontrivial
element in the kernel
of this composition would correspond to a  non-trivial  $L^2$-eigenfunction on
$Y$
with both leading exponents vanishing. Since this   is impossible we have shown
the claim.

We conclude
\begin{prop}
\begin{eqnarray*}
H^p(\Gamma,H^{\lambda}_{-\omega,\Lambda})&=&0, \quad p\ge 2\ ,\\
H^0(\Gamma,H^{\lambda}_{-\omega,\Lambda})&=&\ker(S^{reg}_{-\lambda})\ ,\\
H^1(\Gamma,H^{\lambda}_{-\omega,\Lambda})&=&\coker(S^{reg}_{-\lambda})\ .
\end{eqnarray*}
Moreover,
$\dim H^0(\Gamma,H^{\lambda}_{-\omega,\Lambda}) = \dim
H^1(\Gamma,H^{\lambda}_{-\omega,\Lambda})=\dim\ker\:
(S^{reg}_{-\lambda})_{|V}<\infty$.
\end{prop}

\section{The case $\lambda=-1/2,-3/2,\dots$}\label{uio}
Let $t$ be the parameter of $S^1$ and $d:=d/dt$.
Recall the exact sequences
$$0\rightarrow F_k \rightarrow
H^{k/2}\stackrel{d^k}{\rightarrow}H^{-k/2}\rightarrow F_k\rightarrow 0\ ,$$
$k\in 2\nat_0+1$, of $G$-modules, where $F_k$ is the finite-dimensional
representation of $G$ of dimension $k$.
Define the $\Gamma$-modules $M_k$, $N_k$ by
\begin{eqnarray*}
&&0\rightarrow
H^{k/2}_{-\omega,\Lambda}\stackrel{d^k}{\rightarrow}H^{
-k/2}_{-\omega,\Lambda}\rightarrow M_k\rightarrow 0\\
&&0\rightarrow N_k\rightarrow
H^{k/2}_{-\omega}(\Omega)\stackrel{d^k}{\rightarrow}
H^{-k/2}_{-\omega}(\Omega)\rightarrow 0\ ,
\end{eqnarray*}
where $H^\lambda_{-\omega}(\Omega)$ denotes the hyperfunctions on $\Omega$ with
the $\Gamma$-module
structure given by $\lambda$.
The short exact sequence of complexes
$$\begin{array}{ccccccccc}
&&0&&0&&0&&\\
&&\downarrow&&\downarrow&&\downarrow&&\\
0&\rightarrow&H^{k/2}_{-\omega,\Lambda}&\rightarrow&H^{k/2}_{
-\omega}&\rightarrow&H^{k/2}_{-\omega}(\Omega)&\rightarrow&
0                              \\
&&d^k\downarrow&&d^k\downarrow&&d^k\downarrow&&\\
0&\rightarrow&H^{-k/2}_{-\omega,
\Lambda}&\rightarrow&H^{-k/2}_{-\omega}&\rightarrow&H^{-k/2}_{
-\omega}(\Omega)&\rightarrow&0\\
&&\downarrow&&\downarrow&&\downarrow&&\\
&&0&&0&&0&&\end{array}$$
induces the exact sequence
\begin{equation}\label{haha}0\rightarrow F_k\rightarrow N_k\rightarrow M_k
\rightarrow F_k\rightarrow 0\ .\end{equation}
Since $H^{k/2}_{-\omega,\Lambda}$ has trivial $\Gamma$-cohomology by
Proposition
\ref{erst}, we have
$$H^\ast(\Gamma,M_k)\cong H^\ast(\Gamma, H^{-k/2}_{-\omega,\Lambda})$$
leading to the problem of computing the cohomology of $M_k$.

We assume that $\Omega\not=\emptyset$ and that $\Gamma$ is non-trivial.
Then $\Gamma$ is a free group
$$\Gamma=\langle a_1,\dots,a_g,b_1,\dots,b_g,\sigma_1,\dots,\sigma_{t-1}\rangle
\ ,$$
where $g$ is the genus of the surface $Y$ and $t$ is the number of boundary
components,
i.e., $\dim H^0(B)$. The generators $\sigma_i$ correspond to the boundary
circles.
The missing circle $\sigma_t$ can be expressed in terms of the other
generators.
Let $\Sigma_i\cong \Z$ be the group generated by $\sigma_i$.
\begin{lem}
As a representation of $\Gamma$ we have
$N_k=\oplus_{i=1}^t \Ind_{\Sigma_i}^\Gamma F_{k|\Sigma_i}$.
Moreover, the embedding $F_k\rightarrow N_k$ is given by $u:F_k \ni
v\mapsto\oplus_{i=1}^t (\gamma\mapsto \gamma^{-1} v)$.
\end{lem}
\proof
As a vector space, $N_k=\oplus_{\mbox{\scriptsize components of }\Omega} F_k$,
where we identify $F_k$ with the kernel of $d^k$, i.e., with the polynomials of
order at most $k-1$. The $\Gamma$-action is easy to check.
\hB
The cohomology of $N_k$ can be computed using the Shapiro Lemma \cite{brown82}:
$$H^\ast(\Gamma,N_k)=H^\ast(\Gamma,\oplus_{i=1}^t \Ind_{\Sigma_i}^\Gamma
F_{k|\Sigma_i})=\oplus_{i=1}^t H^\ast(\Sigma_i,F_{k|\Sigma_i})\ .$$
Now $\sigma_i$ is hyperbolic and thus has a unique fixed line in $F_k$.
It follows that $H^p(\Sigma_i,F_{k|\Sigma_i})=\C$ for $p=0,1$ and
$H^p(\Sigma_i,F_{k|\Sigma_i})=0$ for $p\ge 2$.
We  finally obtain
$$H^\ast(\Gamma,N_k)=H^\ast(B)\ .$$
{}From the exact sequence (\ref{haha}) we see that $H^\ast(\Gamma,M_k)$
is finite-dimensional and $\chi(\Gamma,M_k)=\chi(\Gamma,N_k)=0$.
We now split the sequence (\ref{haha}) into two short exact sequences
\begin{eqnarray*}
&&0\rightarrow F_k\rightarrow M_k\rightarrow R_k\rightarrow 0\ ,\\
&&0\rightarrow R_k\rightarrow N_k\rightarrow F_k\rightarrow 0\ ,
\end{eqnarray*}
where $R_k$ is a certain $\Gamma$-module.
We obtain two long exact sequences
\begin{eqnarray}
&&0\rightarrow H^0(\Gamma,F_k)\rightarrow H^0(B)\rightarrow
H^0(\Gamma,R_k)\rightarrow\nonumber\\
&&\rightarrow H^1(\Gamma,F_k)\rightarrow H^1(B)\rightarrow
H^1(\Gamma,R_k)\rightarrow 0\label{ser1}\\[0.5cm]
&&
0\rightarrow H^0(\Gamma,R_k) \rightarrow H^0(\Gamma,M_k)\rightarrow
H^0(\Gamma,F_k)\rightarrow\nonumber\\
&&\rightarrow H^1(\Gamma,R_k)\rightarrow H^1(\Gamma,M_k)\rightarrow
H^1(\Gamma,F_k)\rightarrow 0\ .\label{ser2}
\end{eqnarray}
In the following discussion we will assume that $\Gamma$ is non-abelian.
Then for $k\ge 2$ we have
$$H^0(\Gamma,F_k)=0\quad \mbox{  and  }\quad \dim H^1(\Gamma,F_k)=(2g-2+t)k\
,$$
while in case $k=1$ we have $$\dim H^0(\Gamma,F_1)=1  \mbox{  and  }  \dim
H^1(\Gamma,F_1)=2g+t-1\ .$$
Let $h^1:=\dim H^1(\Gamma,M_k)$ and $q$ be the dimension of the image of
$u_\ast : H^1(\Gamma,F_k)\rightarrow H^1(B)$.
Then we can read off from (\ref{ser1}) that
$\dim H^1(\Gamma,R_k)= t-q$.
In case $k\ge 2$ from (\ref{ser2}) we obtain
\begin{equation}\label{ee3}
h^1=t-q+(2g-2+t)k \ . \end{equation}
In case $k=1$ we claim
\begin{equation}\label{tr1} h^1=t-q+2g+t-1\ .\end{equation}
We must show that
$H^0(\Gamma,M_1)\rightarrow H^0(\Gamma,F_1)$ is surjective.
This amounts to construct an invariant  hyperfunction one-form $\omega$ on
$S^1$ with non-trivial integral and support in $\Lambda$.
In fact, $H^{-1/2}_{-\omega}$ can be identified with the hyperfunction
one-forms
and the map $H^{-1/2}_{-\omega}\rightarrow F_1=\C$ is the integral of the
one-form over $S^1$.
Let $\chi$ be a locally constant function on $B$. Lift it to a locally constant
function $\tilde{\chi}$ on $\Omega$.
Let $\R^1\rightarrow S^1$ be the universal cover and $\hat{\chi}$ be the lift
of $\tilde{\chi}$ to $\R^1$.
Assume that $0\in\R^1$ projects to a point of  $\Lambda$.
Consider $[0,1]$ as a fundamental domain of that cover. We define the function
$\chi_1(t):=\hat{\chi}(t)+k$ for all $t\in\R^1$, $t-k\in (0,1)$, projecting  to
a point of $\Omega$.
Then $d\chi_1$ is a distribution on $\R^1$ carried by the lift of the limit
set.
It projects down to $S^1$ to give the desired $\Gamma$-invariant hyperfunction
(in fact distribution) one-form with integral one.
Thus (\ref{tr1}) holds for $k=1$.

In order to compute $h^1$ it remains to compute $t-q$.
\begin{lem}
We have  $t-q=1$ for $k=1$ or $g=0$
and $t-q=0$ for $k\ge 2$, $g\not=0$.
\end{lem}
 \proof
Let $V$ be a $\Gamma$-module.
Recall the space of group cochains
$$C^i(\Gamma,V)=\{\phi:
\underbrace{\Gamma\times\dots\times\Gamma}_{i\times}\rightarrow V\}$$
and the boundary operators $\partial^0:C^0(\Gamma,V)\rightarrow C^1(\Gamma,V)$,
$(\partial v)(g)=gv-v$, $v\in C^0(\Gamma,V)=V$ and
$\partial^1:C^1(\Gamma,V)\rightarrow C^2(\Gamma,V)$,
$(\partial^1\phi)(g_1,g_2)=g_1(\phi(g_2)-\phi(g_1^{-1}))-\phi(g_1g_2)$.
 A one-cocycle is uniquely determined by its values
on the generators. Since $\Gamma$ is free, the values on the generators can be
prescribed arbitrarily.
Hence
$$Z^1(\Gamma,V)=\C^{2g+t-1}\otimes V\ .$$

Let $\phi\in Z^1(\Gamma,F_k)$ be given by its values on the generators
$\phi(a_i),\phi(b_i),\phi(\sigma_i)$.
Then $u \phi \in Z^1(\Gamma,N_k)$ is given by
$(u\phi)(g)=\oplus_{i=1}^t(h\mapsto h^{-1}\phi(g))$.
This cocycle is a boundary iff the following equations have a solution
$\oplus_{i=1}^t\psi_i\in N_k$:
\begin{eqnarray}
h^{-1}\phi(g)&=&\psi_i(g^{-1}h)-\psi_i(h),\quad i=1,\dots,t,\quad\forall
g,h\in\Gamma\label{er1}\\
(1-\sigma_i^{-1})\psi_i(e)&=&\phi(\sigma_i),\quad i=1,\dots,t\label{er2}\ .
\end{eqnarray}
Here the first equation (\ref{er1}) encodes the boundary map
while the second (\ref{er2}) ensures that $\oplus_{i=1}^t\psi_i\in N_k$.
Since $\phi$ was a cocycle there are always functions $\psi_i:\Gamma\rightarrow
F_k$
solving (\ref{er1}) and the solutions are determined by $\psi_i(e)$.
Thus $u\phi$ is a boundary iff the equations
\begin{equation}\label{et3}(1-\sigma_i^{-1})v=\phi(\sigma_i),\quad i=1,\dots,t\
,\end{equation}
are solvable.
Let us first discuss the cases $k=1$ or $g=0$.
$F_1$ is the trivial representation and (\ref{et3}) is not solvable if
$\phi(\sigma_i)\not=0$.
We claim that for $k=1$ or $g=0$ the value
$\phi(\sigma_t)$ depends on $\phi(\sigma_i)$, $i=1,\dots,t-1$.

The case $g=0$ is obvious. Consider $g\ge 1$.
Then $\sigma_t^{-1}=[a_1,b_1]\dots[a_g,b_g]\sigma_1\dots\sigma_{t-1}$.
A group one-cocycle $\phi$ in the trivial representation vanishes on
commutators.
Thus $\phi([a_1,b_1]\dots[a_g,b_g])=0$. The value
of $\phi(\sigma_1\dots\sigma_{t-1})$ only depends on $\phi(\sigma_i)$,
$i=1,\dots,t-1$.
Thus $\phi(\sigma_t)$ depends on  $\phi(\sigma_i)$, $i=1,\dots,t-1$, too.
This proves the claim.
Since the $\sigma_i$ are hyperbolic we have
$$\dim (\coker(1-\sigma_i^{-1}):F_k\rightarrow F_k)  =1\ .$$
By choosing non-zero $\phi(\sigma_i)$, $i=1,\dots,t-1$, appropriately
we can produce a $t-1$-dimensional subspace in $H^1(\Gamma,N_1)$.
Thus the claim implies $t-q=1$.

We now discuss the case $k\ge 2$, $g\not=0$.
By chosing non-zero $\phi(\sigma_i)$, $i=1,\dots,t-1$, as above
we can produce a $t-1$-dimensional subspace in $H^1(\Gamma,N_k)$.
In order to find a one-dimensional complement, we set $\phi(\sigma_i)=0$
for $i=1,\dots,t-1$, $\phi(a_i)=\phi(b_i)=0$, $i=2,\dots,g$, $\phi(b_1)=0$
and choose $\phi(a_1)$ appropriately.
In fact by the cocycle  equation
\begin{eqnarray*}
\phi(\sigma_t^{-1})&=&[a_1,b_1]\phi([a_2,b_2]\dots[a_g,b_g]
\sigma_1\dots\sigma_{t-1})+\phi([a_1,b_1])\\
&=&\phi([a_1,b_1])\\
&=&(1-a_1b_1a_1^{-1})\phi(a_1)\ .
\end{eqnarray*}
Since $[a_1b_1a_1^{-1},\sigma_t]\not=0$ we can choose $\phi(a_1)$  such that
$(1-a_1b_1a_1^{-1})\phi(a_1)$ represents a non-trivial element in the cokernel
of $(1-\sigma_t)$. Then the equation (\ref{et3}) has no solution for $i=t$.
It follows that $q=t$.
\hB
\begin{prop}
Let $\Gamma$ be non-abelian. Then we have
$$
\dim H^0(\Gamma,H^{-k/2}_{-\omega,\Lambda})= \dim
H^1(\Gamma,H^{-k/2}_{-\omega,\Lambda})=\left\{\begin{array}{cc}(2g-2+t)k, &k\ge
2, g\not=0\\

                                 (2g-2+t)k+1,&k\ge 2,g=0\\

                                 2g +t     ,&k=1
\end{array}\right.\ . $$
\end{prop}
If $\Gamma$ is abelian, then the limit set consists of two points.
$H^\ast(\Gamma,H^\lambda_{-\omega,\Lambda})$ is represented by
derivatives of delta distributions located at $\Lambda$.
The following result is an easy exercise.
\begin{kor} If $\Gamma$ is non-trivial and abelian, then
$$\dim H^0(\Gamma,H^{-k/2}_{-\omega,\Lambda})= \dim
H^1(\Gamma,H^{-k/2}_{-\omega,\Lambda})=2\ .$$
\end{kor}

\section{The case $\Ree(\lambda)=0$, $\lambda\not=0$}

Let $\lambda=\imath\mu\not=0$ and $A_{\imath\mu}:=\Delta-1/4-\mu^2$.
By the theorem of Helgason the Poisson transform provides an embedding
$P_{\imath\mu} :H^{\imath\mu}_{-\omega,\Lambda}\hookrightarrow \cE$.
In order to characterize the range of $P_{\imath\mu}$ we introduce the
following semi-norms:
$q_{W,k}(f):=\|A_{\imath\mu}^kf\|_{B(W)}$, $p_W(f):=\|f\|_{B^\ast(W)}$,
$s_W(f)=p_W(D_\mu f)$,
and $r_{K,i}(f):=\|A_{\imath\mu}^if\|_{L^2(K)}$. The indices run over the
domains
$k\ge 1$, $i\ge 0$. $K$ runs over the compact subsets of $X$ and $W$ over the
sets of the form
$W= M \times (0,a) \subset \R^2_+$, $a>0$, for compact
$M\subset \Omega$, using the coordinates $(x,y)$, $y>0$, of the upper
half-plane model.
The operator $D_\mu$ is given in these coordinates by
$D_\mu :=y\frac{d}{dy}-(\imath \mu+1/2)$. The $B$- and $B^\ast$-norms are
defined by
$\|f\|_{B(W)}:=\sum_{j\ge 0}2^{j/2}\|f\|_{L^2(W\cap \Omega_j)}$,
$\|f\|_{B^\ast(W)}:=\sup_{j\ge 0} 2^{-j/2}\|f\|_{L^2(W\cap \Omega_j)}$,
$\Omega_j:=\{(x,y) \in \R^2_+| -\ln(y)\in[2^{j-1},2^j)\}$, $j\ge 1$, and
$\Omega_0:=\{(x,y)\in \R^2_+|y\ge \ee\}$.
By $B(W)$ and  $B^\ast(W)$ we denote the Banach spaces of functions on $W$ with
finite $B$- or $B^\ast$-norm.
We also introduce the closed subspace ${}^\circ B^\ast(W)=\{f\in
B^\ast(W)|\lim\sup_R\frac{1}{R}\|f\|^2_{L^2(M\times(a,\ee^{-R}))}=0\}$.
These norms are natural in the framework of scattering theory
\cite{hoermander83}, Ch.14.

Let $\cE^\mu_\Lambda$ be the Frechet space of all $f\in\cE$ such that for all
$W$ and $K$
described above
$q_{W,k}(f) <\infty$, $\forall k\ge 1$, $p_W(f)<\infty$, $s_W(f)<\infty$,
$r_{K,i}(f)<\infty$, $\forall i\ge 0$ and
 $D_\mu f \in {}^\circ B^\ast(W)$.
Set  $\cE_\Lambda^\mu(A_{\imath\mu})=\cE^\mu_\Lambda\cap\ker(A_{\imath\mu})$.
\begin{lem}
$$ P_{ \imath\mu}(H^{\imath\mu}_{-\omega,\Lambda}) =
\cE_\Lambda^\mu(A_{\imath\mu})$$
\end{lem}
\proof
Let $P_{ \imath\mu}((x,y),b)=(\frac{y}{y^2+(x-b)^2})^{\imath\mu+1/2}$ be the
kernel defining the Poisson transform.
Then $A_{\imath\mu} P_{ \imath\mu}(.,b)=0$, $P_{ \imath\mu}(.,b)\in B^\ast(W)$
and $D_\mu P_{\imath\mu}(.,b)\in {}^\circ B^\ast(W)$, $W= M \times (0,a)$, for
$b\not\in M$.
Hence $b\to P_{ \imath\mu}(.,b)$ defines an analytic function from a
neighbourhood
of $\Lambda$ to $\cE^\mu_{\Lambda|W}$ (with the obvious definition) for all
$W$.
It follows
$ P_{ \imath\mu}(H^{\imath\mu}_{-\omega,\Lambda})\subset
\cE_\Lambda^\mu(A_{\imath\mu})$.

 We now show that $f\in \cE_\Lambda^\mu(A_{\imath\mu})$  implies
that the boundary value $\beta_{ \imath\mu}(f)$ vanishes on $\Omega\cap\bar{W}$
and hence
$\beta_{ \imath\mu}(\cE_\Lambda^\mu(A))\subset
H^{\imath\mu}_{-\omega,\Lambda}$.
Here $\beta_{ \imath\mu}$ is the boundary value corresponding to the asymptotic
$y^{-\imath\mu+1/2}$.

The argument is similar to the corresponding argument in the proof of Lemma
\ref{volley}. We  again employ the coordinates of the disk model.
Let $\psi\in H^{-\imath\mu}_\omega$. Then
$$\langle \beta_{\imath\mu}(f),\psi\rangle=\lim_{R\to
\infty}\frac{1}{R}\int_0^R f_\psi(t) \ee^{t(1/2-\imath\mu)} dt\ ,$$
where
$$f_\psi(t)=\int_{S^1} f(t,\alpha) \psi(\ee^{\imath\alpha})
\frac{d\alpha}{2\pi}\ ,$$
and the $t$-coordinate is the hyperbolic distance from the origin.
Let $(a,b)\subset \Omega$. Using $(\frac{d}{dt}+\imath\mu+1/2)f=-D_\mu f +
o(t^{-1})$
and $D_\mu f\in {}^\circ B^*$ along $\Omega$ we obtain
\begin{eqnarray*}
\lefteqn{\frac{1}{R}\int_0^R \int_{(a,b)}(\frac{d}{dt}+\imath\mu+1/2)
f(t,\alpha)\ee^{ t(\mp \imath\mu-1/2)} \psi(\ee^{\imath\alpha})
\frac{d\alpha}{2\pi}\ee^t  dt}\hspace{2cm}\\
&\le& (R^{-1/2}\|D_\mu f\|_{L^2((0,R)\times (a,b))} + o(R^{-3/2})  )
\|\psi\|_{L^2((a,b))}\\
& \stackrel{R\to \infty}{\rightarrow} &0\ .
\end{eqnarray*}
By partial integration this implies  on the one hand
$$\lim_{R\to \infty}\frac{1}{R} \ee^{R(1/2+\imath\mu)} \int_{(a,b)} f(R,\alpha)
\psi(\ee^{\imath\alpha}) \frac{d\alpha}{2\pi}=0 $$
and on the other hand
\begin{eqnarray*}
\lefteqn{2\imath \mu \lim_{R\to \infty}\frac{1}{R}\int_0^R\int_{(a,b)}
f(t,\alpha) \psi(\ee^{\imath\alpha}) \frac{d\alpha}{2\pi}
\ee^{t(1/2-\imath\mu)}dt}\hspace{2cm}   \\
&=& -\lim_{R\to \infty}
\frac{1}{R} \ee^{R(1/2-\imath\mu)} \int_{(a,b)} f(R,\alpha)
\psi(\ee^{\imath\alpha}) \frac{d\alpha}{2\pi}   \ .
\end{eqnarray*}
 Thus
$$\lim_{R\to \infty}\frac{1}{R}\int_0^R\int_{(a,b)} f(t,\alpha)
\psi(\ee^{\imath\alpha}) \frac{d\alpha}{2\pi} \ee^{t(1/2-\imath\mu)}dt=0\ .$$
We see that $\beta_{\imath\mu}(f)$ defines a continuous functional on the germs
of analytic sections of $T^{-\imath\mu-1/2}$ on $S^1\setminus (a,b)$.
The argument can now be completed as in Lemma \ref{volley}.
\hB
We introduce another  Frechet space of functions
$$\cE_\Lambda:=\{f\in C^\infty(X)\:|\: q_{W,k}(f)<\infty,\forall W, k\ge 0,
r_{K,i}(f)<\infty,\forall K, i\ge 0\}\ .$$
The next lemma is proved in analogy with \cite{bunkeolbrich947}, Lemma 2.4.
\begin{lem}
$\cE^\mu_\Lambda$, $\cE_\Lambda$ are $\Gamma$-acyclic.
\end{lem}
The following Lemma implies that
\begin{equation}\label{zu}
0\longrightarrow H^{\imath\mu}_{-\omega,\Lambda} \stackrel{P}{\longrightarrow}
\cE^\mu_\Lambda\stackrel{A_{\imath\mu}}{\longrightarrow}\cE_\Lambda
\longrightarrow 0
\end{equation}
is a $\Gamma$-acyclic resolution of $H^{\imath\mu}_{-\omega,\Lambda}$.
\begin{lem}\label{surt}
$A_{\imath\mu}:\cE_\Lambda^\mu\rightarrow \cE_\Lambda$ is surjective.
\end{lem}
\proof
We consider the adjoint operator ${}^tA_{\imath\mu}:
\cE_\Lambda^\prime\rightarrow (\cE_\Lambda^\mu)^\prime$.
We must show that ${}^tA_{\imath\mu}$ is injective and has closed range.
Since $C_c^\infty(X)$ is dense in $\cE_\Lambda$, we can embed
$ \cE_\Lambda^\prime$ into the distributions on $X$.
In fact if $f\in  \cE_\Lambda^\prime$, then it is a distribution
with support in a finite union of fundamental domains of  $\Gamma$
(see the proof of Lemma \ref{wolf} for a similar argument).
If ${}^tA_{\imath\mu} f=0$, then $f$ is real analytic. Since it vanishes on a
non-empty open subset of $X$,
it vanishes identically. Thus ${}^tA_{\imath\mu}$ is injective.

We must show that ${}^tA_{\imath\mu}$ has  closed range.
As in Lemma \ref{wolf} we can restrict the consideration to bounded sets
$B^\prime \subset(\cE_\Lambda^\mu)^\prime$.
There is a subset $U\subset X$ being the union of a compact set $K\subset X$
and $W=(0,b)\times M$, $M$ compact, $M\subset \Omega$, such that $X\setminus U$
is connected and $\supp(h)\subset U$
for all $h\in B^\prime$.
Let ${}^tA_{\imath\mu} f_i=:h_i\in B^\prime$ such that
$h_i\to h$ in $(\cE_\Lambda^\mu)^\prime$. We have $\supp(f_i)\subset U$.
We must find  a $f\in \cE_\Lambda^\prime$ with ${}^tA_{\imath\mu} f=h$.

Let $\cE^\mu_\emptyset,\cE_\emptyset$ be defined like
$\cE_\Lambda^\mu,\cE_\Lambda$
but for the trivial limit set.
Then $f_i\in \cE_\emptyset^\prime$ and ${}^tA_{\imath\mu} f_i\to h$ holds in
$(\cE^\mu_\emptyset)^\prime$. Using a result of Perry \cite{perry87}
(essentially a corollary to Lemma \ref{uni})  we conclude that
$A_{\imath\mu}:\cE^\mu_\emptyset \rightarrow \cE_\emptyset$
is a topological isomorphism.
Thus $f_i\to f={}^tA_{\imath\mu}^{-1}h$ in $(\cE_\emptyset)^\prime$.
Since $\supp(f_i)\subset U$ we also have $\supp(f)\subset U$ and hence $f\in
\cE_\Lambda^\prime$.
Thus the   range of ${}^tA_{\imath\mu}:\cE_\Lambda^\prime\rightarrow
(\cE_\Lambda^\mu)^\prime $
is closed.
 This finishes the proof of surjectivity of $A_{\imath\mu}$ and of the lemma.
\hB
We now consider the $\Gamma$-invariants in the resolution (\ref{zu}).
Using the coordinates (\ref{coor}) of a collar neighbourhood $E$ of infinity
$T_Y:(0,\epsilon)\times B\rightarrow E\subset Y$
we can define the spaces $B(Y)$, $B^\ast(E)$ and ${}^\circ B^\ast(E)$
as in \cite{perry87}.
Let $A_{{\imath\mu},Y}=\Delta_Y-1/4- \mu^2$, where $\Delta_Y$ is the Laplacian
on $Y$
and $D_\mu=\frac{d}{dr}-\imath\mu$.
Then ${}^\Gamma\cE_\Lambda^\mu=\cE_Y^\mu$,${}^\Gamma \cE_\Lambda=\cE_Y$
with
\begin{eqnarray*}
\cE_Y^\mu&:=&\{f\in C^\infty(Y)| \|A_{{\imath\mu},Y}^kf\|_{B(Y)}<\infty,k\ge
1,\| f\|_{B^\ast(E)}<\infty , D_\mu f\in{}^\circ B^\ast(E)\}\\
\cE_Y&:=&\{f\in C^\infty(Y)| \|A_{{\imath\mu},Y}^kf\|_{B(Y)}<\infty,k\ge 0\}
\end{eqnarray*}

It again follows from the results of \cite{perry87}  that
$A_{{\imath\mu},Y}:\cE_Y^\mu\rightarrow\cE_Y$
is a  topological isomorphism.  Thus
$\ker(A_{{\imath\mu},Y})=\coker(A_{{\imath\mu},Y})=0$.
\begin{prop}\label{imaga}
For $0\not=\lambda$, $\Ree(\lambda)=0$, we have
$$H^\ast(\Gamma,H^\lambda_{-\omega,\Lambda})=0\ .$$
\end{prop}

\section{The scattering matrix near $\Ree(\lambda)=0$}\label{scats}

In Section \ref{sss} we constructed a continuation of the scattering matrix
up to the imaginary axis.
In this section we first show that it is continuous at the imaginary line.
Then we employ the functional equation in order to provide
the meromorphic continuation to all of $\C$.
Consider the operator $A_{\lambda,Y}=\Delta_Y-1/4+\lambda^2$.
For $\Ree(\lambda)>0$ it is invertible on $L^2(Y)$.
If $\lambda$ approaches the imaginary axis we consider its inverse
on the slightly smaller space $B(Y)$.
\begin{lem}[Perry, \cite{perry87}, Prop. 5.5]\label{uni}
The inverse $A_{\lambda,Y}^{-1}:B(Y)\rightarrow  B^\ast(Y)$
is a weakly continuous family of bounded operators
uniformly bounded on compact subsets of
$\lambda\in\C\setminus(\{0\}\cup \Sigma)$, $\Ree(\lambda)\ge 0$.
\end{lem}
We now employ the constructions and notations introduced in Section \ref{sss}.
Let $\phi\in C^\infty(B)={}^\Gamma H^{1/2}_\infty(\Omega)$.
We can consider $\phi_\lambda:=\varho^{\lambda-1/2}\phi\in {}^\Gamma
H^\lambda_\infty(\Omega)$.
We want to construct a hyperfunction extension $\Phi_\lambda\in
H^\lambda_{-\omega}$.
The exact sequence
$$0\longrightarrow H^\lambda_{-\omega,\Lambda}\longrightarrow
H^\lambda_{-\omega,(\Lambda) }\stackrel{res}{\longrightarrow}
H^\lambda_{\infty}(\Omega)\longrightarrow 0$$
implies the long exact sequence
\begin{eqnarray*}
&&0\longrightarrow {}^\Gamma H^\lambda_{-\omega,\Lambda}\longrightarrow
{}^\Gamma H^\lambda_{-\omega,(\Lambda) }\stackrel{res}{\longrightarrow}
{}^\Gamma H^\lambda_{\infty}(\Omega)\longrightarrow\\
&&\longrightarrow H^1(\Gamma,H^\lambda_{-\omega,\Lambda})\longrightarrow\dots\
{}.
\end{eqnarray*}
{}From Proposition \ref{imaga} it follows that $\Phi_\lambda$ exists for
$\Ree(\lambda)\ge 0$, $\lambda\not\in \Sigma$, $\lambda\not=0$.
As we have seen in Section \ref{sss} $\Phi_\lambda$ is holomorphic near the
imaginary axis.
We show that $\Phi_\lambda$ is continuous up to $\Ree(\lambda)=0$,
$\lambda\not=0$.
We construct the family of approximate eigenfunctions $f_\lambda$ and note
that $g_\lambda:=A_{\lambda,Y}f_\lambda\in B(Y)$.
Lemma \ref{uni} shows that $F_\lambda=f_\lambda-A_{\lambda,Y}^{-1}g_\lambda$
extends continuously to $\Ree(\lambda)=0$, $\lambda\not= 0$, as a family of
eigenfunctions in $C^\infty(Y)$ (weakly continuous in $B^\ast(Y)$).
Lifting to the universal cover and taking the boundary value
we see that the family of hyperfunctions $\Phi_\lambda$ is continuous up to the
imaginary axis, too.
We define the right limit $S_{\lambda+0}$ of the scattering matrix for
$\lambda\not=0$, $\Ree(\lambda)=0$,
by
$$S_{\lambda+0}\phi_\lambda:=(res\circ J_\lambda )\Phi_{\lambda+0}\ .$$
Then by definition
$S(\lambda)=\varho^{1/2+\lambda}S_\lambda\varho^{\lambda-1/2}$
is continuous in the strong topology as an operator on $C^\infty(B)$
for $\Ree(\lambda)\ge 0$ and $\lambda\not=0$, $\lambda\not\in \Sigma$ and
$\lambda\not\in \nat_0$.

The Knapp-Stein intertwining operators satisfy the functional equation
$$
J_\lambda\circ J_{-\lambda} = - \frac{\ctg(\pi\lambda)}{2\lambda}\ .
$$
It follows for $\Ree(\lambda)=0$, $\lambda\not=0$,
\begin{eqnarray*}
S_{\lambda+0}\circ S_{-\lambda+0}\phi_\lambda &=&\lim_{\epsilon\to 0} (res\circ
J_{\lambda+\epsilon})(res\circ J_{-\lambda+\epsilon}\Phi_{-\lambda+\epsilon}
)_{\mbox{\scriptsize extended}}\\
&=&res\circ  J_\lambda  \circ J_{-\lambda}\Phi_{\lambda+0}\\
&=&-\frac{\ctg(\pi \lambda)}{2\lambda}\phi_\lambda\ .
\end{eqnarray*}
We define $S_\lambda$ for $\Ree(\lambda)< 0$ by
$$S_\lambda:=-\frac{\ctg(\pi \lambda)}{2 \lambda}S_{-\lambda}^{-1}\ .$$
A standard application of the Fredholm theory for Fr\'echet spaces
\cite{grothendieck56} similar as in \cite{patterson76} shows that $S_\lambda$
is meromorphic with at most finite-dimensional singularities at
$\lambda\not\in\frac{1}{2}\nat_0$.
It remains  to show continuity of $S_\lambda$ at $\Ree(\lambda)=0$.
In fact we have
$$
S_{\lambda+0}-S_{\lambda- 0}=
S_{\lambda+0}+ \frac{\ctg(\pi \lambda)}{2\lambda}S_{-\lambda+0}^{-1}=0\ .
$$
Thus we have shown
\begin{lem}\label{sctti}
$S(\lambda)$ is a meromorphic family of elliptic
pseudodifferential operators on $C^\infty(B)$.
It has at most finite-dimensional poles at the spectral points
$\lambda\in\Sigma$ and
in the set of resonances $\Ree(\lambda)<0$  with
$H^0(\Gamma,H^{ \lambda}_{-\omega,\Lambda})\not=0$.
There are further singularities at $\lambda\in \nat_0$.
\end{lem}
\begin{kor}
The set of $\lambda\in\C$ with
$H^\ast(\Gamma,H^\lambda_{-\omega,\Lambda})\not=0$
is discrete.
\end{kor}
We now can find an analytic continuation of $\Phi_\lambda$.
For $\Ree(\lambda)<0$, $\lambda$ not a resonance or in $\frac{1}{2}\Z$,
define
$$\Phi_\lambda:=-2\lambda\tg(\pi\lambda) J_{-\lambda}\circ res^{-1}\circ
S_\lambda \phi_\lambda\ .$$
Then for $\Ree(\lambda)=0$, $\lambda\not=0$, we have
\begin{eqnarray*}
\Phi_{\lambda-0}&=&-2\lambda\tg(\pi\lambda)J_{-\lambda+0}\circ res^{-1}\circ
S_{\lambda-0} \phi_{\lambda-0}\\
&=&-2\lambda\tg(\pi\lambda)J_{-\lambda}\circ J_{\lambda}\circ res^{-1}
\phi_{\lambda}\\
&=& res^{-1}\phi_\lambda\\
&=&\Phi_{\lambda+0}\ .
\end{eqnarray*}
It follows that  $\Phi_\lambda$ is a meromorphic family of hyperfunctions.
We obtain the Eisenstein series by applying the Poisson transform:
$E(\phi,\lambda)=P_\lambda \Phi_\lambda$.
\begin{prop}
The Eisenstein series has a meromorphic continuation to all of $\C\setminus
\{0\}$
with poles in $\Sigma$, the set of resonances, and the negative half-integers.
The Eisenstein series satisfies the functional equation
$$E(\phi,\lambda)=
\sqrt{2}\frac{\Gamma(1/2-\lambda)}{\Gamma(-\lambda)}E(S_\lambda\phi,-\lambda)\
.$$
\end{prop}
The functional equation follows easily from the corresponding functional
equation of the Poisson transform
\cite{kashiwarakowataminemuraokamotooshimatanaka78}.
The meromorphic continuation of the scattering matrix and the
Eisenstein series for surfaces was first
obtained by Patterson
\cite{patterson75},\cite{patterson761},\cite{patterson761}.
Generalizations to  higher dimensions can be found in
\cite{patterson89}, \cite{perry89}, \cite{mandouvalos88}, \cite{mandouvalos86}.
Since a generalization with
more technical details can be found in \cite{bunkeolbrich953} in the
present paper we kept the arguments concerning the scattering matrix and the
Eisenstein series sketchy.

\section{Summary}
Let $\Gamma$ be a torsion-free Fuchsian group of the second kind without
parabolic elements. By $H^\lambda_{-\omega,\Lambda}$
we denote the hyperfunction vectors of the principal series
representation  $H^\lambda_{-\omega}$ with support in the limit set
of $\Gamma$.
For $\lambda\not=0$ the cohomology $H^\ast(\Gamma,
H^\lambda_{-\omega,\Lambda})$
is finite-dimensional and has vanishing Euler characteristic.
In greater detail:
\begin{itemize}
\item For $\Ree(\lambda)>0$, we have
$$H^0(\Gamma,H^\lambda_{-\omega,\Lambda})=H^1(\Gamma,
H^\lambda_{-\omega,\Lambda})=  \ker_{L^2}(\Delta_Y-1/4+\lambda^2)\ ,$$
where $\Delta_Y$ is the Laplace-Beltrami operator on the hyperbolic surface
$Y:=\Gamma\backslash H^2$. Let $\Sigma$ be the set of $\lambda$,
$\Ree(\lambda)>0$, with $H^0(\Gamma,H^\lambda_{-\omega,\Lambda})$ non-trivial.
\item For $\Ree(\lambda)=0$, $\lambda\not=0$, we have
$$H^\ast(\Gamma,H^\lambda_{-\omega,\Lambda})=0\ .$$
\item The case $\lambda=0$ is still unknown.
\item For $\Ree(\lambda)<0$, $\lambda\not\in -\frac12\nat$,
$-\lambda\not\in\Sigma$,
we have
 \begin{eqnarray}
H^0(\Gamma,H^\lambda_{-\omega,\Lambda})&=&\ker(S_{-\lambda})\label{summ}\\
H^1(\Gamma,H^\lambda_{-\omega,\Lambda})&=&\coker(S_{-\lambda})\ ,\nonumber
\end{eqnarray}
where $S_\lambda$ is the scattering matrix. We have
$\dim\: H^0(\Gamma,H^\lambda_{-\omega,\Lambda})
=\dim\:H^1(\Gamma,H^\lambda_{-\omega,\Lambda})$.
\item For $\Ree(\lambda)<0$, $-\lambda \in\Sigma$, (\ref{summ}) holds with
$S_{-\lambda}$
replaced by $S^{reg}_{-\lambda} $ defined
  in Section \ref{s4}.
\item For $ \Ree(\lambda)<0$, $-\lambda \in  \nat$,  (\ref{summ})  holds with
$S_{-\lambda}$ replaced by the renormalized $\tilde{S}_{-\lambda}$.
\item For $\lambda=-k/2$, $k$ odd, we have for non-abelian $\Gamma$
$$
\dim H^0(\Gamma,H^{-k/2}_{-\omega,\Lambda})= \dim
H^1(\Gamma,H^{-k/2}_{-\omega,\Lambda})=\left\{\begin{array}{cc}(2g-2+t)k, &k\ge
2, g\not=0\\

                                 (2g-2+t)k+1,&k\ge 2,g=0\\

                                 2g +t     ,&k=1
\end{array}\right. $$
and for abelian $\Gamma$
$$\dim H^0(\Gamma,H^{-k/2}_{-\omega,\Lambda})= \dim
H^1(\Gamma,H^{-k/2}_{-\omega,\Lambda})=2 \ .$$
Here $g$ is the genus and $t$ is the number of boundary components of $Y$.
 \end{itemize}

\bibliographystyle{plain}


\end{document}